\documentclass[12pt,a4paper]{amsart} 
\usepackage[utf8]{inputenc}
\usepackage[english]{babel}

\usepackage{amsmath}
\usepackage{amsfonts}
\usepackage{amssymb}
\usepackage{color}
\usepackage{ifthen}

\usepackage{amsthm}
\newtheorem{theorem}{Theorem}
\newtheorem{proposition}{Proposition}

\usepackage{graphics}
\usepackage{epsfig}

\usepackage{url}
\usepackage{hyperref}

\usepackage[a4paper,top=2.5cm,bottom=2.8cm,
left=3.2cm,right=3.2cm]{geometry}
\markleft{\hfill \textsc{Inequalities For Distances Between Triangle Centers} \hfill}
\long\def\void#1{}

\begin{document}
\bigskip
\medskip

\begin{center}
	{\Large \textbf{Inequalities For Distances}} \\
	\medskip
	{\Large \textbf{Between Triangle Centers}} \\
	\bigskip

	\textsc{Stanley Rabinowitz} \\

	545 Elm St Unit 1,  Milford, New Hampshire 03055, USA \\
	e-mail: \href{mailto:stan.rabinowitz@comcast.net}{stan.rabinowitz@comcast.net}\\
	web: \url{http://www.StanleyRabinowitz.com/} \\
	
\end{center}
\bigskip

\textbf{Abstract.} In his seminal paper on triangle centers, Clark Kimberling
made a number of conjectures concerning the distances between triangle centers.
For example, if $D(i,j)$ denotes the distance between triangle centers $X_i$ and $X_j$,
Kimberling conjectured that $D(6,1)\leq D(6,3)$ for all triangles.
We use symbolic mathematics techniques to prove these conjectures.
In addition, we prove stronger results, using best-possible
constants, such as $D(6,1)\leq (2-\sqrt3)D(6,3)$.

\medskip
\textbf{Keywords.} triangle geometry, triangle centers, inequalities, computer-discovered mathematics, Blundon's Fundamental Inequality GeometricExplorer.

\medskip
\textbf{Mathematics Subject Classification (2020).} 51M04, 51-08.


\section{Introduction}
\label{section:intro}

Let $X_n$ denote the $n$th named triangle center as cataloged
in the Encyclopedia of Triangle Centers \cite{ETC}.
Let $X_iX_j$ denote the distance between $X_i$ and $X_j$.
We will also write this as $D(i,j)$.

In his seminal paper on triangle centers \cite{Kimberling}, Clark Kimberling
made a number of conjectures concerning the distances between pairs of triangle centers.
For example,
Kimberling conjectured that $D(6,1)\leq D(6,3)$ for all triangles.

He also conjectured the truth of many chains of inequalities, such as the following.
$$X_3X_9\leq X_3X_{10}\leq X_3X_2\leq X_3X_{12}\leq X_3X_7\leq X_3X_4.$$

Kimberling reached these conjectures by using a computer to examine 10,740
different shaped triangles and
numerically computing the coordinates for the centers.
Upon determining that the inequality held for each of these 10,740 triangles,
he then conjectured that the inequality was true for all triangles.

With the advances in computers and symbolic algebra systems,
it is now possible to prove these conjectures using exact symbolic computation.

\section{Barycentric Coordinates}
\label{section:methodology}

We use barycentric coordinates in this study.
The barycentric coordinates for triangle centers $X_1$ through $X_{20}$
in terms of the sides of the triangle, $a$, $b$, and $c$, are shown in Table \ref{table:centers},
where
$$S=\frac12 \sqrt{(a + b - c) (a - b + c) (-a + b + c) (a + b + c)}.$$
Only the first barycentric coordinate is given, because if $f(a,b,c)$
is the first barycentric coordinate for a point $P$, then the barycentric coordinates
for $P$ are
$$\Bigl(f(a,b,c):f(b,c,a):f(c,a,b)\Bigr).$$
These were derived from \cite{ETC}.

\begin{table}[h!t]
\caption{Barycentric coordinates for the first 20 centers}
\label{table:centers}
\begin{center}
\begin{tabular}{|l|l|}
\hline
n&first barycentric coordinate for $X_n$\\
\hline
1&$a$\\ \hline
2&1\\ \hline
3&$a^2(a^2-b^2-c^2)$\\ \hline
4&$(a^2 + b^2 - c^2) (a^2 - b^2 + c^2)$\\ \hline
5&$c^4-a^2 b^2 + b^4 - a^2 c^2 - 2 b^2 c^2$\\ \hline
6&$a^2$\\ \hline
7&$(a + b - c) (a - b + c)$\\ \hline
8&$a - b - c$\\ \hline
9&$a (a - b - c)$\\ \hline
10&$b + c$\\ \hline
11&$(b - c)^2 (-a + b + c)$\\ \hline
12&$(a + b - c) (a - b + c) (b + c)^2$\\ \hline
13&$a^4 - 2 (b^2 - c^2)^2 + a^2 (b^2 + c^2 + 2 \sqrt3 S)$\\ \hline
14&$a^4 - 2 (b^2 - c^2)^2 + a^2 (b^2 + c^2 - 2 \sqrt3 S)$\\ \hline
15&$a^2 (\sqrt3 (a^2 - b^2 - c^2) - 2 S)$\\ \hline
16&$a^2 (\sqrt3 (a^2 - b^2 - c^2) + 2 S)$\\ \hline
17&$(a^2+b^2-c^2+2 \sqrt3 S) (a^2-b^2+c^2+2 \sqrt3 S)$\\ \hline
18&$(a^2 + b^2 - c^2 - 2 \sqrt3 S) (a^2 - b^2 + c^2 - 2 \sqrt3 S)$\\ \hline
19&$a (a^2 + b^2 - c^2) (a^2 - b^2 + c^2)$\\ \hline
20&$3 a^4 - 2 a^2 b^2 - b^4 - 2 a^2 c^2 + 2 b^2 c^2 - c^4$\\ \hline
\end{tabular}
\end{center}
\end{table}

To find the distance between two centers, we used the following formula
which comes from \cite{Grozdev}.

\begin{proposition}
Given two points $P = (u_1, v_1, w_1)$ and $Q = (u_2, v_2, w_2)$ in normalized barycentric coordinates. Denote $x = u_1 -u_2$, $y = v_1 -v_2$ and $z = w_1 -w_2$.
Then the distance between $P$ and $Q$ is
$$\sqrt{-a^2yz-b^zx-c^2xy}.$$
\end{proposition}

\newpage
\section{Graphs}
\label{section:graphs}

For $n$, $i$, and $j$ ranging from 1 to 20, we used Algorithm B from \cite{Rabinowitz} to check every inequality of the form
$D(n,i)\leq D(n,j)$. Algorithm B is based on Blundon's Fundamental Inequality \cite{Blundon}.
Figure $n$ shows a graph of the results.
An arrow from node $i$ to node $j$ means that $D(n,i)\leq D(n,j)$ for all triangles.
No arrow means the inequality does not hold for all triangles.
Since we used exact symbolic computations, these results are theorems and not conjectures.
To avoid radicals, we replaced inequalities of the form 
$D(a,b)\leq D(c,d)$ by the equivalent inequality
$D(a,b)^2\leq D(c,d)^2$.

\begin{figure}[h!t]
\centering
\includegraphics{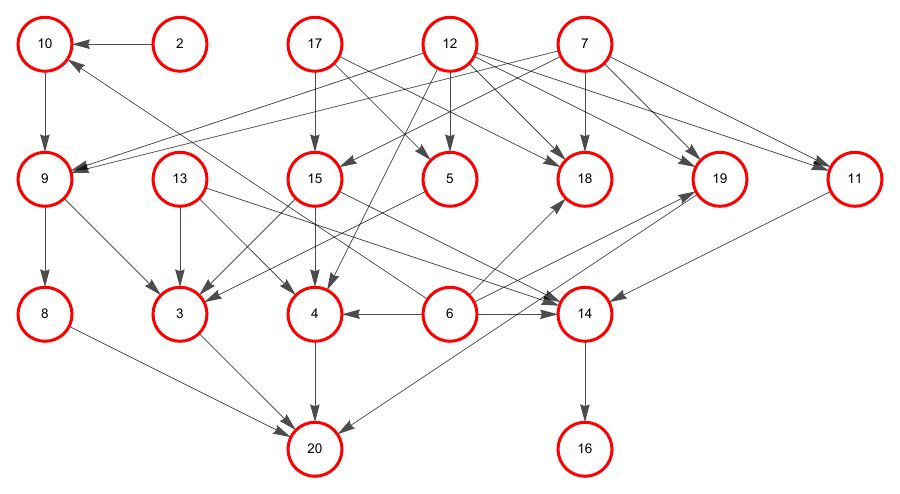}
\caption{$X_1$ inequalities. An arrow from $i$ to $j$ means $X_1X_i\leq X_1X_j$.}
\label{fig:X1}
\end{figure}

\begin{figure}[h!t]
\centering
\includegraphics{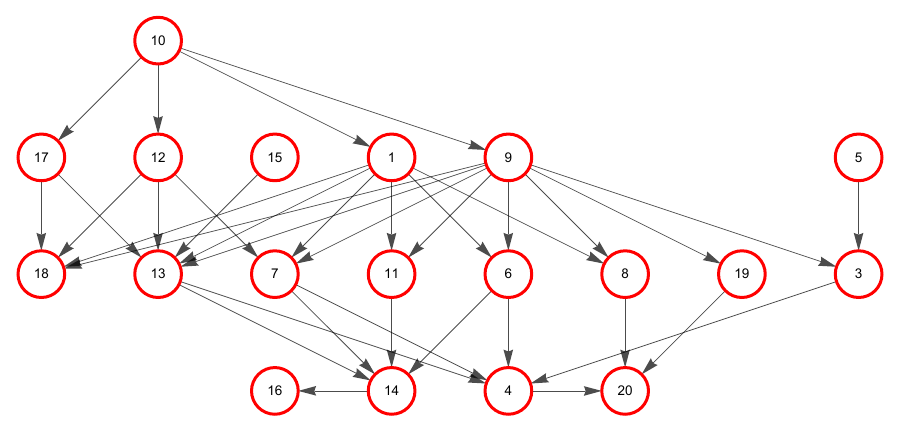}
\caption{$X_2$ inequalities. An arrow from $i$ to $j$ means $X_2X_i\leq X_2X_j$.}
\label{fig:X2}
\end{figure}

\newpage

\begin{figure}[h!t]
\centering
\includegraphics[width=0.9\linewidth]{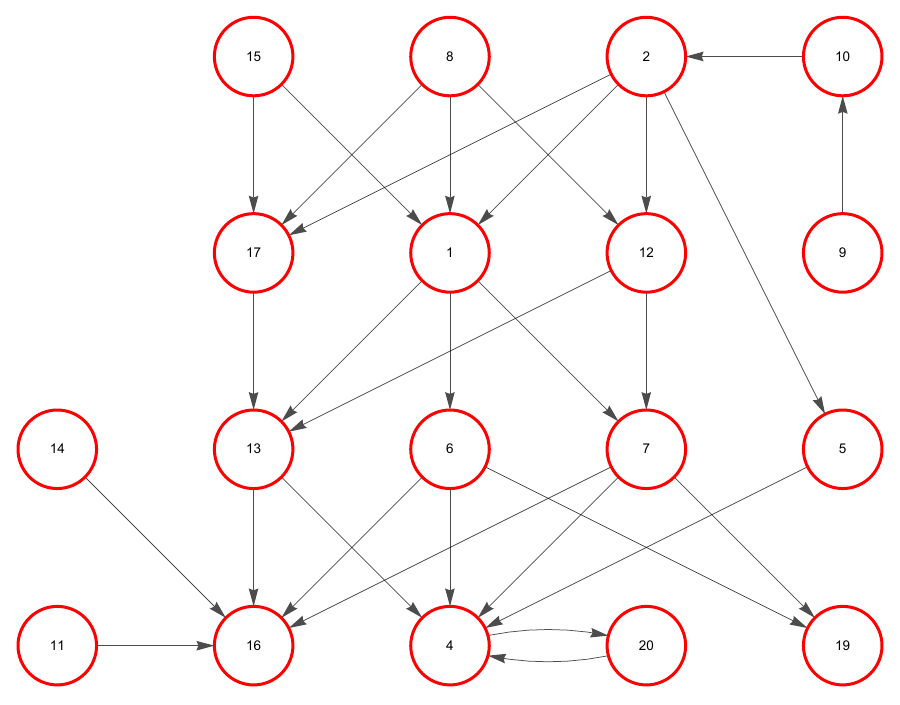}
\caption{$X_3$ inequalities. An arrow from $i$ to $j$ means $X_3X_i\leq X_3X_j$.}
\label{fig:X3}
\end{figure}

\begin{figure}[h!t]
\centering
\includegraphics[width=0.9\linewidth]{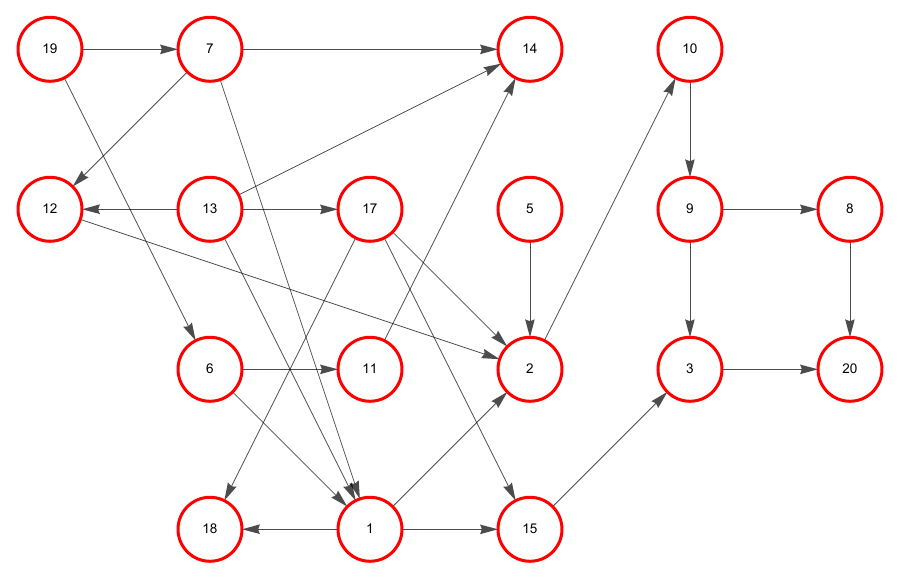}
\caption{$X_4$ inequalities. An arrow from $i$ to $j$ means $X_4X_i\leq X_4X_j$.}
\label{fig:X4}
\end{figure}

\newpage

\begin{figure}[h!t]
\centering
\includegraphics[width=0.8\linewidth]{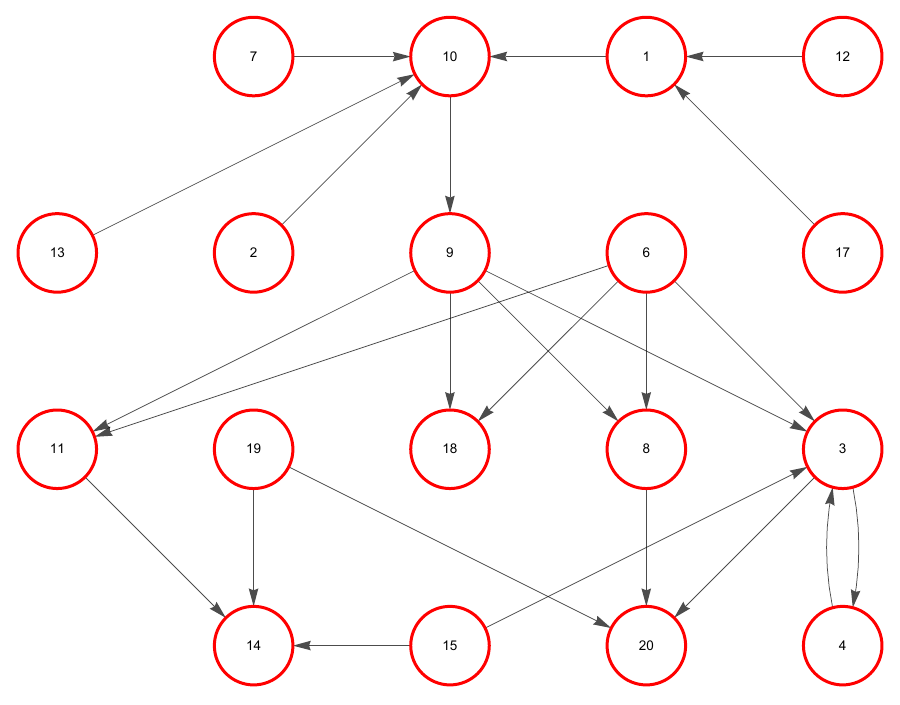}
\caption{$X_5$ inequalities. An arrow from $i$ to $j$ means $X_5X_i\leq X_5X_j$.}
\label{fig:X5}
\end{figure}

\begin{figure}[h!t]
\centering
\includegraphics[width=0.8\linewidth]{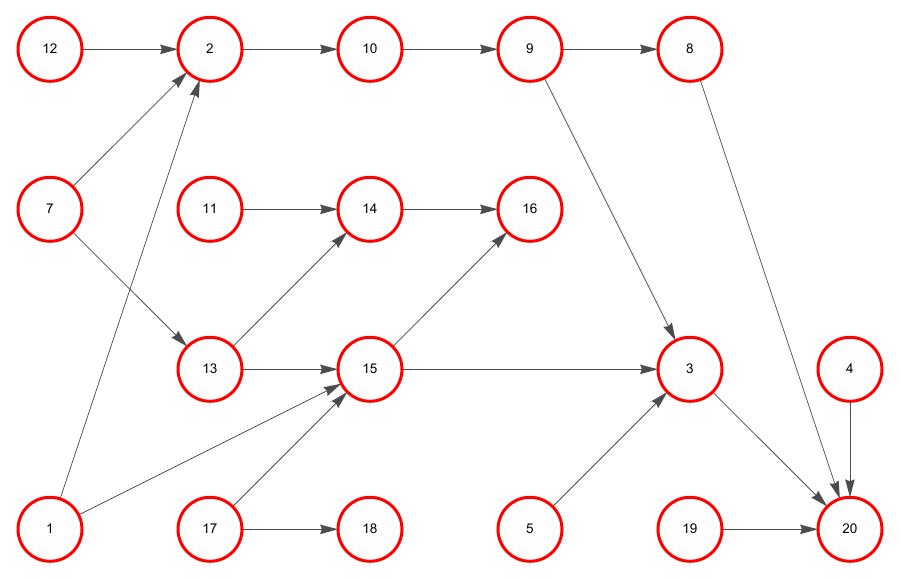}
\caption{$X_6$ inequalities. An arrow from $i$ to $j$ means $X_6X_i\leq X_6X_j$.}
\label{fig:X6}
\end{figure}

\newpage

\begin{figure}[h!t]
\centering
\includegraphics[width=0.7\linewidth]{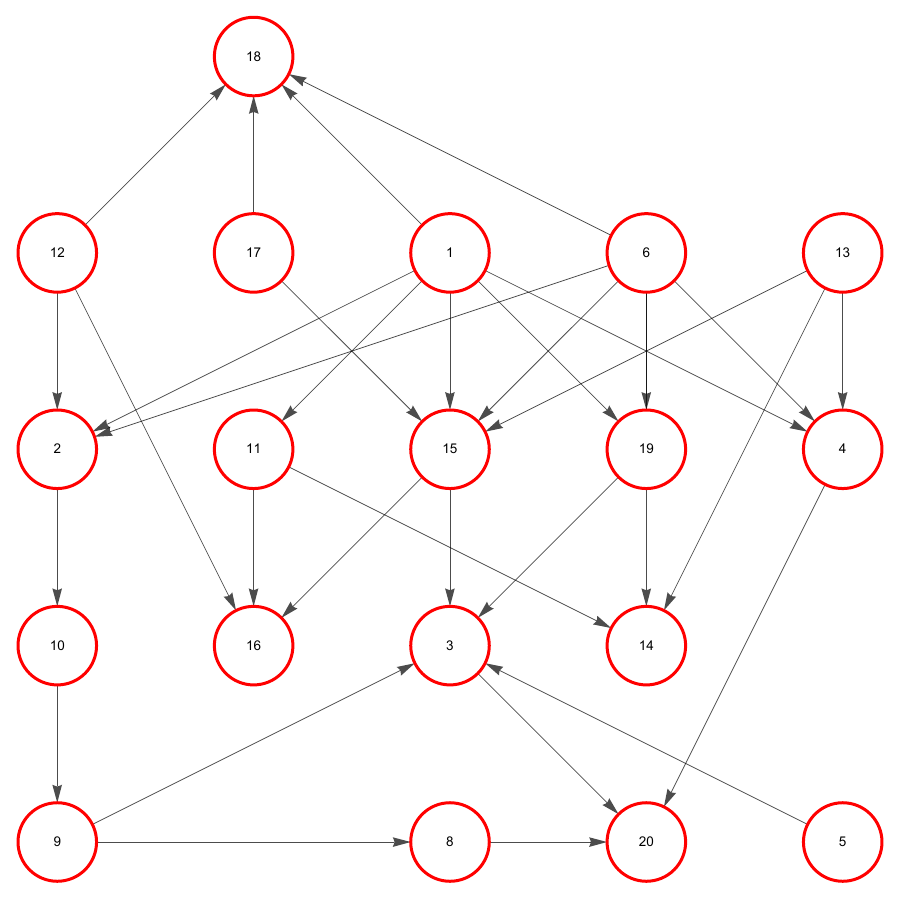}
\caption{$X_7$ inequalities. An arrow from $i$ to $j$ means $X_7X_i\leq X_7X_j$.}
\label{fig:X7}
\end{figure}

\begin{figure}[h!t]
\centering
\includegraphics[width=0.9\linewidth]{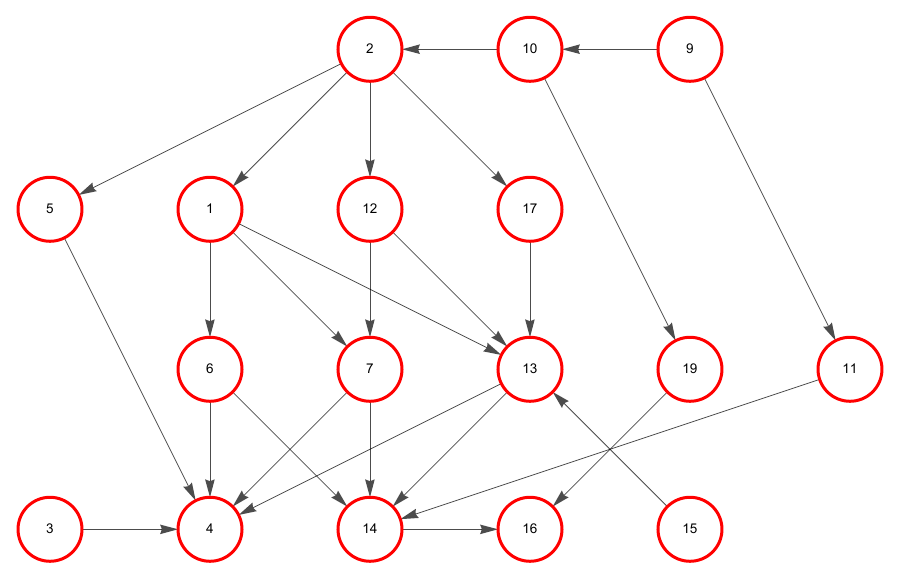}
\caption{$X_8$ inequalities. An arrow from $i$ to $j$ means $X_8X_i\leq X_8X_j$.}
\label{fig:X8}
\end{figure}

\newpage

\begin{figure}[h!t]
\centering
\includegraphics[width=0.8\linewidth]{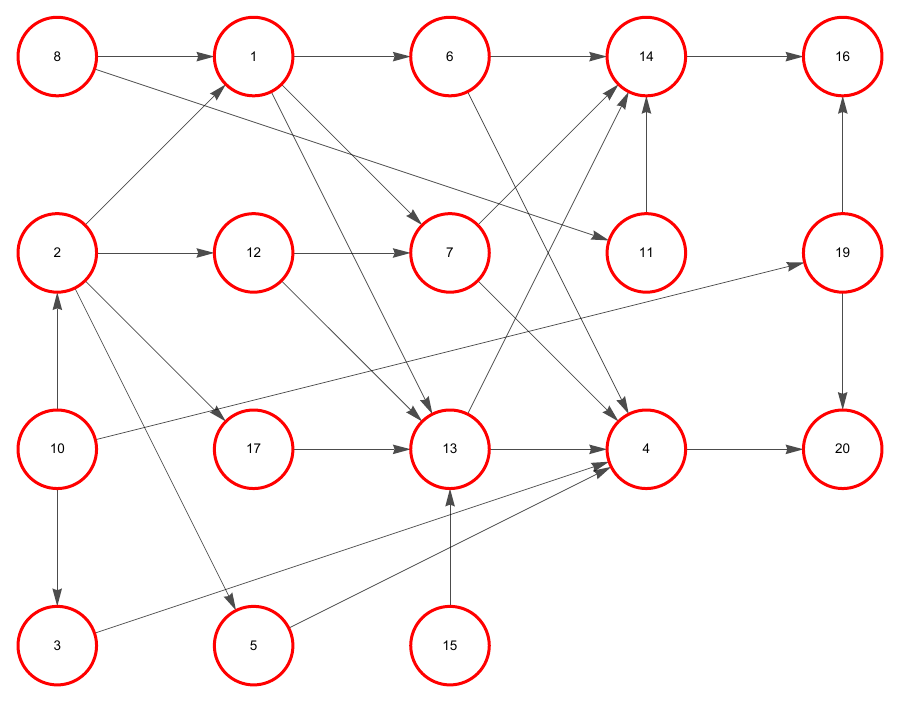}
\caption{$X_9$ inequalities. An arrow from $i$ to $j$ means $X_9X_i\leq X_9X_j$.}
\label{fig:X9}
\end{figure}

\begin{figure}[h!t]
\centering
\includegraphics[width=1\linewidth]{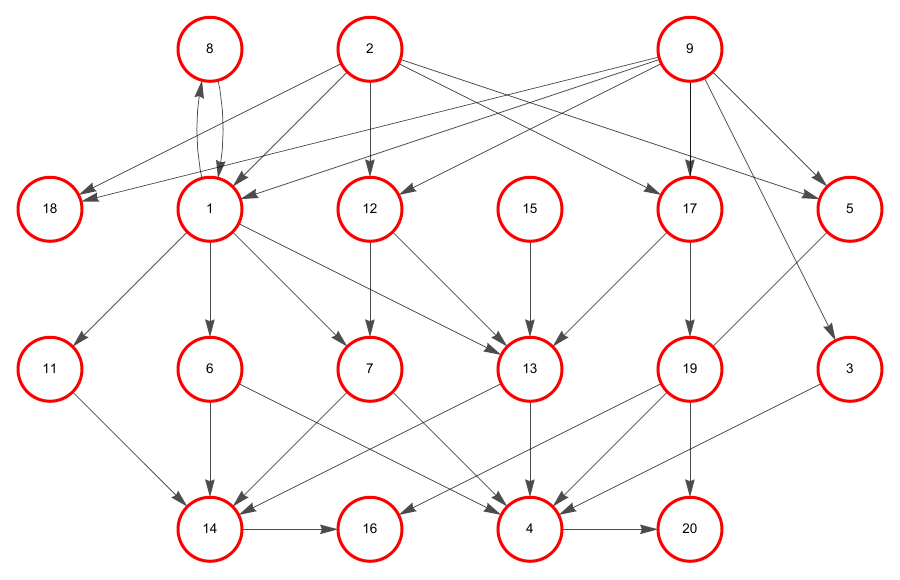}
\caption{$X_{10}$ inequalities. An arrow from $i$ to $j$ means $X_{10}X_i\leq X_{10}X_j$.}
\label{fig:X10}
\end{figure}

\newpage

\begin{figure}[h!t]
\centering
\includegraphics[width=1\linewidth]{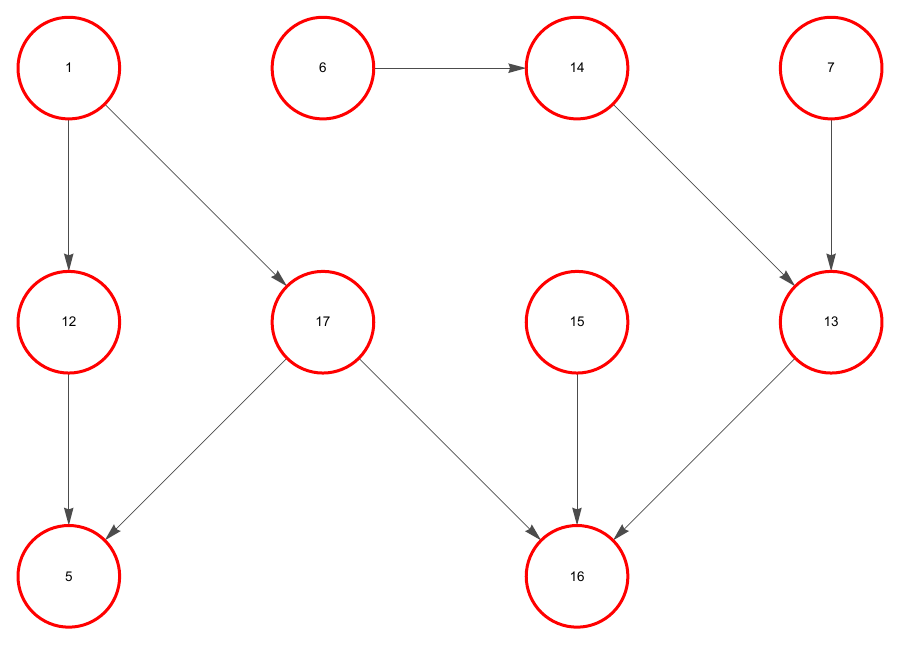}
\caption{$X_{11}$ inequalities. An arrow from $i$ to $j$ means $X_{11}X_i\leq X_{11}X_j$.}
\label{fig:X11}
\end{figure}

\begin{figure}[h!t]
\centering
\includegraphics[width=1\linewidth]{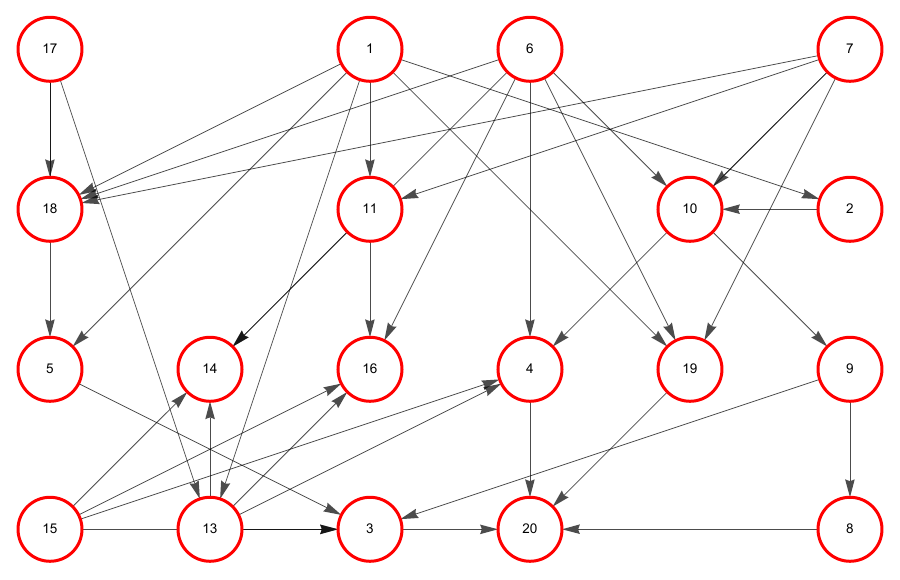}
\caption{$X_{12}$ inequalities. An arrow from $i$ to $j$ means $X_{12}X_i\leq X_{12}X_j$.}
\label{fig:X12}
\end{figure}

\newpage

\begin{figure}[h!t]
\centering
\includegraphics[width=0.7\linewidth]{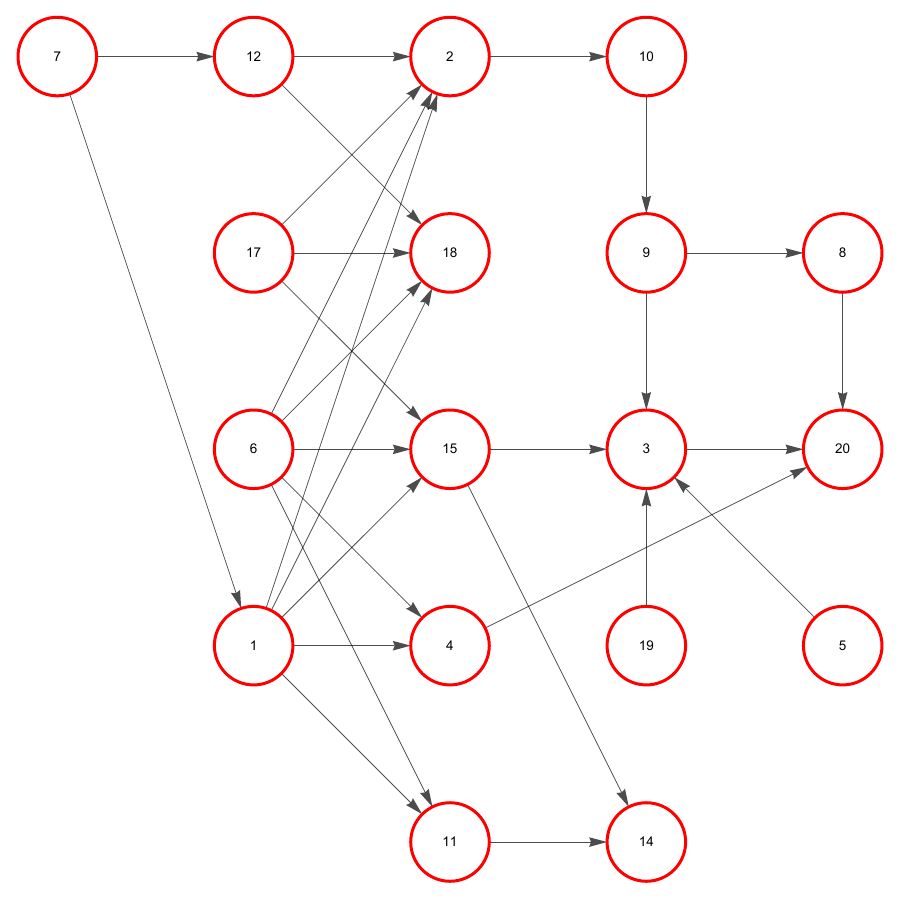}
\caption{$X_{13}$ inequalities. An arrow from $i$ to $j$ means $X_{13}X_i\leq X_{13}X_j$.}
\label{fig:X13}
\end{figure}

\begin{figure}[h!t]
\centering
\includegraphics[width=0.8\linewidth]{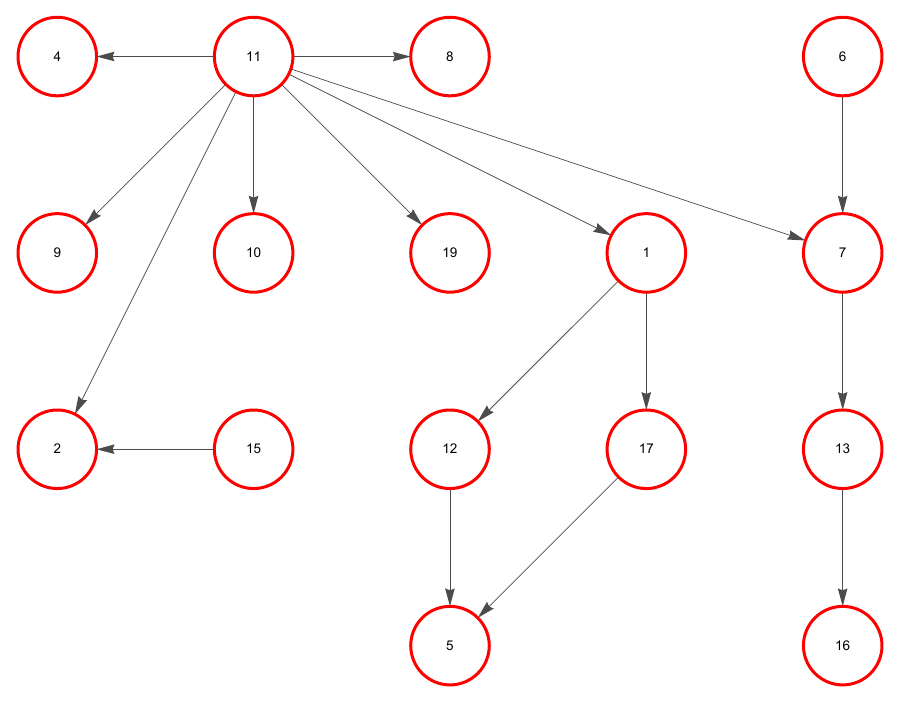}
\caption{$X_{14}$ inequalities. An arrow from $i$ to $j$ means $X_{14}X_i\leq X_{14}X_j$.}
\label{fig:X14}
\end{figure}

\newpage

\begin{figure}[h!t]
\centering
\includegraphics[width=0.9\linewidth]{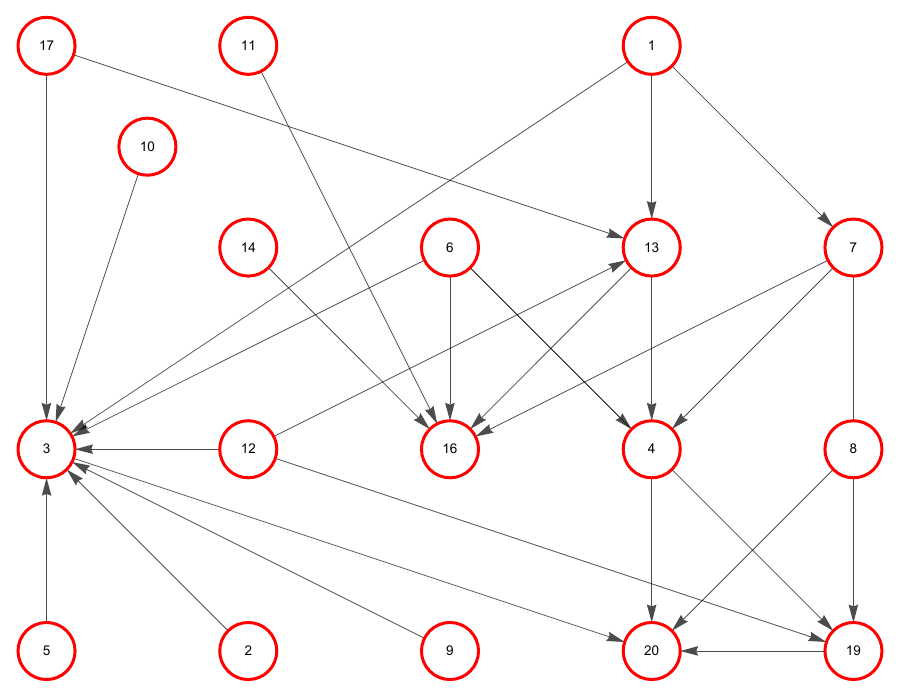}
\caption{$X_{15}$ inequalities. An arrow from $i$ to $j$ means $X_{15}X_i\leq X_{15}X_j$.}
\label{fig:X15}
\end{figure}

\begin{figure}[h!t]
\centering
\includegraphics[width=0.8\linewidth]{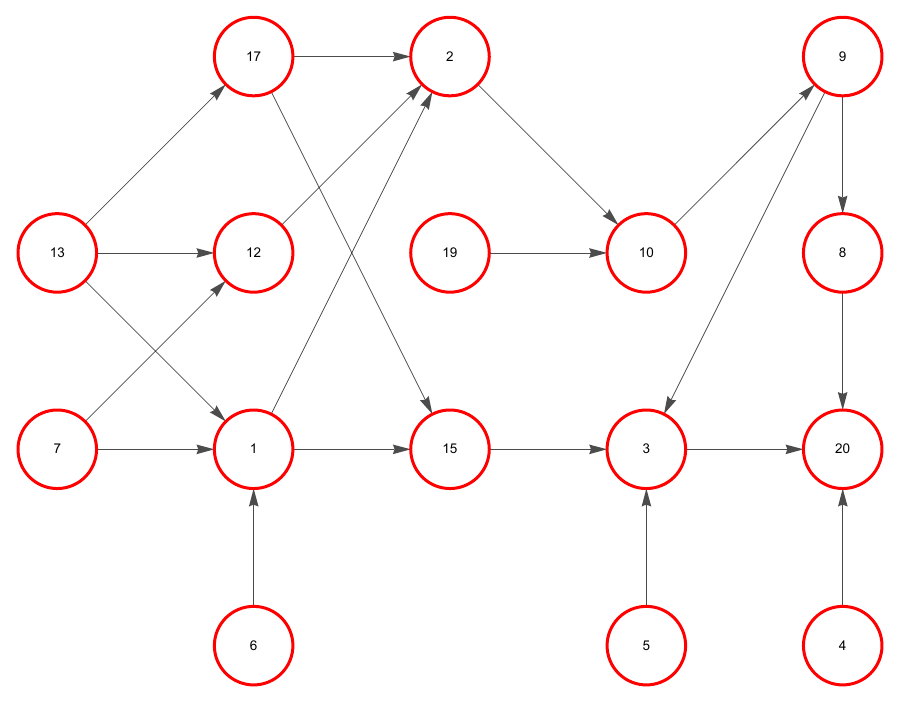}
\caption{$X_{16}$ inequalities. An arrow from $i$ to $j$ means $X_{16}X_i\leq X_{16}X_j$.}
\label{fig:X16}
\end{figure}

\newpage

\begin{figure}[h!t]
\centering
\includegraphics[width=0.8\linewidth]{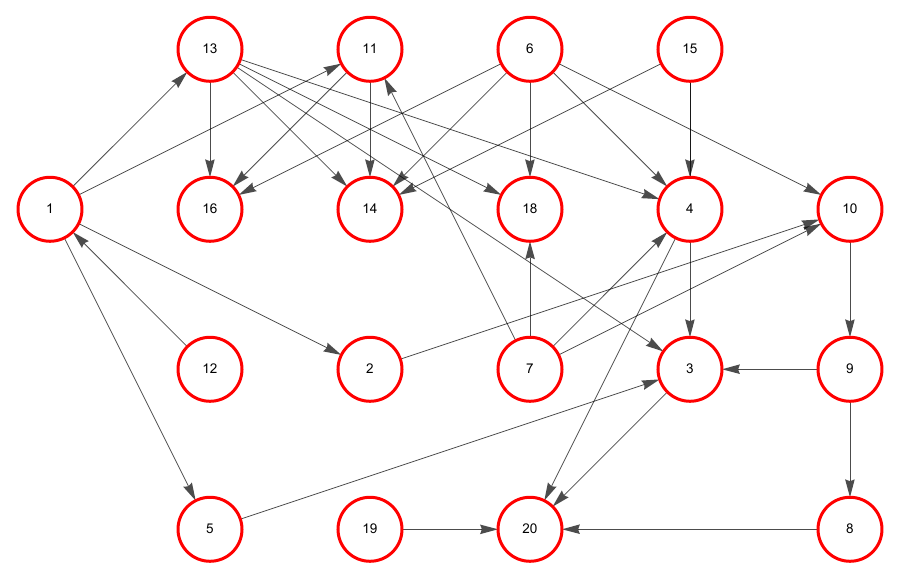}
\caption{$X_{17}$ inequalities. An arrow from $i$ to $j$ means $X_{17}X_i\leq X_{17}X_j$.}
\label{fig:X17}
\end{figure}

There were no inequalities found for $n=18$. In other words, there were no inequalities
of the form
$D(18,i)\leq D(18,j)$
for any $i$ and $j$ with $1\leq i\leq 20$, $1\leq j\leq 20$, $i\neq 18$,
$j\neq 18$, and $i\neq j$.

\begin{figure}[h!t]
\centering
\includegraphics[width=0.08\linewidth]{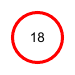}
\caption{There are no inequalities of the form $X_{18}X_i\leq X_{18}X_j$.}
\label{fig:X18}
\end{figure}

\begin{figure}[h!t]
\centering
\includegraphics[width=0.8\linewidth]{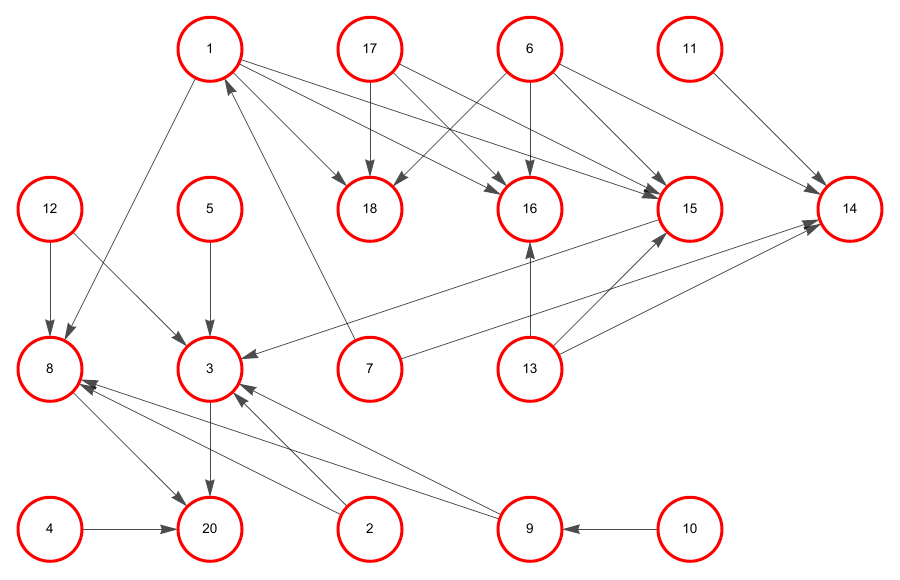}
\caption{$X_{19}$ inequalities. An arrow from $i$ to $j$ means $X_{19}X_i\leq X_{19}X_j$.}
\label{fig:X19}
\end{figure}

\newpage

\begin{figure}[h!t]
\centering
\includegraphics[width=0.9\linewidth]{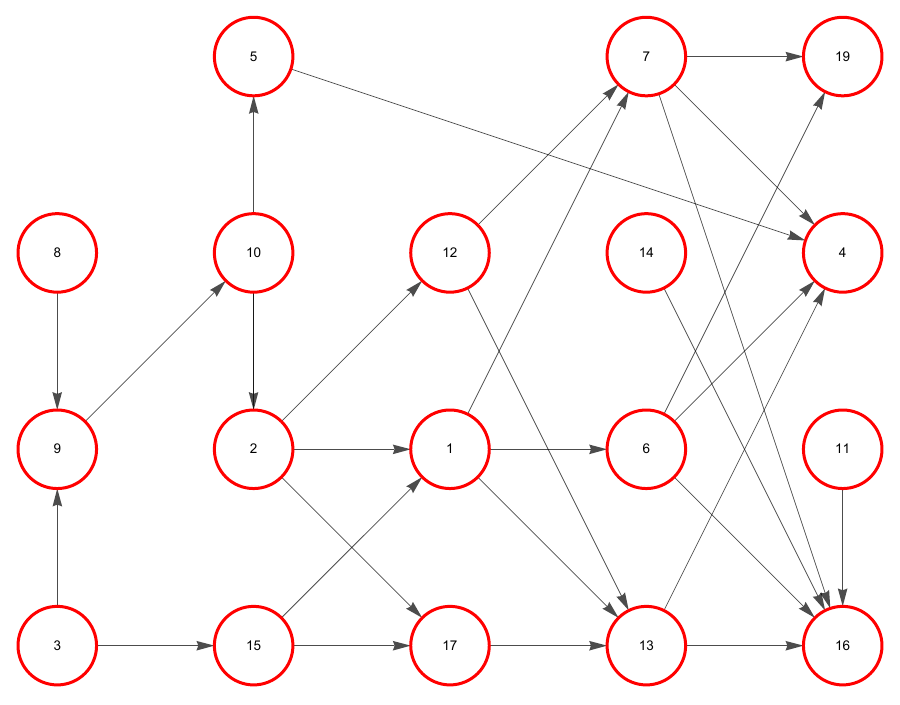}
\caption{$X_{20}$ inequalities. An arrow from $i$ to $j$ means $X_{20}X_i\leq X_{20}X_j$.}
\label{fig:X20}
\end{figure}

Examining these graphs, we note that there are a few loops.
An arrow from $i$ to $j$ and an arrow from $j$ to $i$ in Figure $n$
means that $D(n,i)\leq D(n,j)$ and $D(n,j)\leq D(n,i)$.
This implies that $D(n,i)=D(n,j)$.
Three such equalities were noticed:
$D(1,10)=D(8,10)$, $D(3,4)=D(3,20)$, and $D(3,5)=D(4,5)$.

These equalities were noticed by Kimberling in \cite[Table~5.4]{Kimberling}.
These correspond to the (now) well-known facts that in all triangles,
$X_{10}$ is the midpoint of $\overline{X_1X_8}$, 
$X_{3}$ is the midpoint of $\overline{X_4X_{20}}$, and 
$X_{5}$ is the midpoint of $\overline{X_3X_4}$.

Since we only investigated inequalities between distances formed by three triangle centers,
this does not mean that we can conclude that there aren't any other equalities
of the form $D(i_1,i_2)=D(i_3,i_4)$, where $i_1$, $i_2$, $i_3$, and $i_4$ are
all distinct.

To check for such equalities, we ran a separate Mathematica program that examined
all distances of the form $D(i,j)$ where $i$ and $j$ are distinct integers
between 1 and 20, looking for duplicate distances. No new equalities were found.
This lets us state the following result.

\begin{proposition}
The only pairs of centers from among the first 20 centers that have equal distances
are the following.
$$
\begin{aligned}
D(1,10)&=D(8,10)\\
D(3,4)&=D(3,20)\\
D(3,5)&=D(4,5)
\end{aligned}
$$
\end{proposition}

\newpage
\section{Bounds}
\label{section:bounds}

Some of the inequalities from Section \ref{section:graphs} can be strengthened.
For example, from Figure \ref{fig:X6}, one can see that $D(6,2)\leq D(6,10)$.
However, the stronger inequality
$$D(6,2)\leq \frac13\left(1+\sqrt2\right)D(6,10)$$
is true. To find the best such inequalities,
we applied Algorithm K from \cite{Rabinowitz} to every inequality of the form
$$D(n,i)\leq k D(n,j)\quad\textrm{or}\quad D(n,i)\geq k D(n,j)$$
for $n$, $i$, and $j$ ranging from 1 to 10 with $i<j$
to find the smallest (resp. largest) constant $k$ making the inequality true.
The results are given below, shown as lower and upper bounds for $\displaystyle\frac{D(n,i)}{D(n,j)}$.
Lower bounds of 0 and upper bounds of $\infty$ are omitted.

For example, $\displaystyle 0\leq\frac{D(1,2)}{D(1,4)}\leq \infty$ would mean that Algorithm K
proved that there is no constant $k>0$ such that $\displaystyle k\leq \frac{D(1,2)}{D(1,4)}$ is true for all triangles,
and that there is no constant $k$ such that $\displaystyle \frac{D(1,2)}{D(1,4)}\leq k$ is true for all triangles.

\begin{theorem}
The following bounds involving distances from $X_1$ hold for all triangles.
\end{theorem}

\begin{center}
\begin{tabular}{c|c|c}
\begin{minipage}{0.3\linewidth}
\begin{center}
$\begin{aligned}
 &\frac{D(1,2)}{D(1,3)}\leq \frac{2}{3} \\
 \frac{2 \sqrt{2}}{3}\leq &\frac{D(1,2)}{D(1,6)} \\
 \frac{1}{3}\leq &\frac{D(1,2)}{D(1,7)} \\
 &\frac{D(1,2)}{D(1,8)}=\frac{1}{3} \\
 \frac{1}{3}\leq &\frac{D(1,2)}{D(1,9)}\leq \frac{2}{3} \\
 &\frac{D(1,2)}{D(1,10)}=\frac{2}{3} \\
 \frac{1}{2}\leq &\frac{D(1,3)}{D(1,4)} \\
 2\leq &\frac{D(1,3)}{D(1,5)} \\ \\ \\
\end{aligned}$
\end{center}
\end{minipage}
&
\begin{minipage}{0.26\linewidth}
\begin{center}
$\begin{aligned}
 1+\sqrt{3}\leq &\frac{D(1,3)}{D(1,6)} \\
 \frac{3}{2}+\sqrt{2}\leq &\frac{D(1,3)}{D(1,7)} \\
 \frac{1}{2}\leq &\frac{D(1,3)}{D(1,8)} \\
 1\leq &\frac{D(1,3)}{D(1,9)} \\
 1\leq &\frac{D(1,3)}{D(1,10)} \\
 1\leq &\frac{D(1,4)}{D(1,6)} \\
 2\leq &\frac{D(1,4)}{D(1,7)} \\
 \frac{1}{3}\leq &\frac{D(1,6)}{D(1,7)} \\ \\ \\
\end{aligned}$
\end{center}
\end{minipage}
&
\begin{minipage}{0.3\linewidth}
\ \ \ $\begin{aligned}
 &\frac{D(1,6)}{D(1,8)}\leq \frac{1}{2 \sqrt{2}} \\
 &\frac{D(1,6)}{D(1,9)}\leq \frac{1}{2} \\
 &\frac{D(1,6)}{D(1,10)}\leq \frac{1}{\sqrt{2}} \\
 &\frac{D(1,7)}{D(1,8)}\leq 1 \\
 &\frac{D(1,7)}{D(1,9)}\leq 1 \\
 &\frac{D(1,7)}{D(1,10)}\leq 2 \\
 1\leq &\frac{D(1,8)}{D(1,9)}\leq 2 \\
 &\frac{D(1,8)}{D(1,10)}=2 \\
 1\leq &\frac{D(1,9)}{D(1,10)}\leq 2 \\
\end{aligned}$
\end{minipage}
\end{tabular}
\end{center}

\newpage

\begin{theorem}
The following bounds involving distances from $X_2$ hold for all triangles.
\end{theorem}

\begin{center}
\begin{tabular}{c|c|c}
\begin{minipage}{0.3\linewidth}
\begin{center}
$\begin{aligned}
 &\frac{D(2,1)}{D(2,3)}\leq 2 \\
 &\frac{D(2,1)}{D(2,4)}\leq 1 \\
 &\frac{D(2,1)}{D(2,5)}\leq 4 \\
 6 \sqrt{2}-8\leq &\frac{D(2,1)}{D(2,6)}\leq 1 \\
 \frac{1}{4}\leq &\frac{D(2,1)}{D(2,7)}\leq 1 \\
 &\frac{D(2,1)}{D(2,8)}=\frac{1}{2} \\
 \frac{1}{2}\leq &\frac{D(2,1)}{D(2,9)}\leq 2 \\
 &\frac{D(2,1)}{D(2,10)}=2 \\
 &\frac{D(2,3)}{D(2,4)}=\frac{1}{2} \\
 &\frac{D(2,3)}{D(2,5)}=2 \\
 \frac{1}{2}\leq &\frac{D(2,3)}{D(2,6)} \\
 \frac{1}{2}\leq &\frac{D(2,3)}{D(2,7)} \\
\end{aligned}$
\end{center}
\end{minipage}
&
\begin{minipage}{0.25\linewidth}
\begin{center}
$\begin{aligned}
 \frac{1}{4}\leq &\frac{D(2,3)}{D(2,8)} \\
 1\leq &\frac{D(2,3)}{D(2,9)} \\
 1\leq &\frac{D(2,3)}{D(2,10)} \\
 &\frac{D(2,4)}{D(2,5)}=4 \\
 1\leq &\frac{D(2,4)}{D(2,6)} \\
 1\leq &\frac{D(2,4)}{D(2,7)} \\
 \frac{1}{2}\leq &\frac{D(2,4)}{D(2,8)} \\
 2\leq &\frac{D(2,4)}{D(2,9)} \\
 2\leq &\frac{D(2,4)}{D(2,10)} \\
 \frac{1}{4}\leq &\frac{D(2,5)}{D(2,6)} \\
 \frac{1}{4}\leq &\frac{D(2,5)}{D(2,7)} \\
 \frac{1}{8}\leq &\frac{D(2,5)}{D(2,8)} \\
\end{aligned}$
\end{center}
\end{minipage}
&
\begin{minipage}{0.35\linewidth}
\ \ \ $\begin{aligned}
 \frac{1}{2}\leq &\frac{D(2,5)}{D(2,9)} \\
 \frac{1}{2}\leq &\frac{D(2,5)}{D(2,10)} \\
 \frac{1}{2}\leq &\frac{D(2,6)}{D(2,7)}\leq \frac{3}{2} \\
 \frac{1}{2}\leq &\frac{D(2,6)}{D(2,8)}\leq \frac{4+3 \sqrt{2}}{8} \\
 1\leq &\frac{D(2,6)}{D(2,9)}\leq 3 \\
 2\leq &\frac{D(2,6)}{D(2,10)}\leq 2+\frac{3}{\sqrt{2}} \\
 \frac{1}{2}\leq &\frac{D(2,7)}{D(2,8)}\leq 2 \\
 &\frac{D(2,7)}{D(2,9)}=2 \\
 2\leq &\frac{D(2,7)}{D(2,10)}\leq 8 \\
 1\leq &\frac{D(2,8)}{D(2,9)}\leq 4 \\
 &\frac{D(2,8)}{D(2,10)}=4 \\
 1\leq &\frac{D(2,9)}{D(2,10)}\leq 4 \\
\end{aligned}$
\end{minipage}
\end{tabular}
\end{center}

\newpage

\begin{theorem}
The following bounds involving distances from $X_3$ hold for all triangles.
\end{theorem}

\begin{center}
\begin{tabular}{c|c|c}
\begin{minipage}{0.4\linewidth}
\begin{center}
$\begin{aligned}
 1\leq &\frac{D(3,1)}{D(3,2)}\leq 3 \\
 \frac{1}{3}\leq &\frac{D(3,1)}{D(3,4)}\leq 1 \\
 \frac{2}{3}\leq &\frac{D(3,1)}{D(3,5)}\leq 2 \\
 \sqrt{3}-1\leq &\frac{D(3,1)}{D(3,6)}\leq 1 \\
 \frac{1}{17} \left(7+4 \sqrt{2}\right)\leq &\frac{D(3,1)}{D(3,7)}\leq 1 \\
 1\leq &\frac{D(3,1)}{D(3,8)} \\
 1\leq &\frac{D(3,1)}{D(3,9)} \\
 1\leq &\frac{D(3,1)}{D(3,10)} \\
 &\frac{D(3,2)}{D(3,4)}=\frac{1}{3} \\
 &\frac{D(3,2)}{D(3,5)}=\frac{2}{3} \\
 \frac{1}{3}\leq &\frac{D(3,2)}{D(3,6)}\leq 1 \\ \\ \\
\end{aligned}$
\end{center}
\end{minipage}
&
\begin{minipage}{0.26\linewidth}
\begin{center}
$\begin{aligned}
 \frac{1}{3}\leq &\frac{D(3,2)}{D(3,7)}\leq 1 \\
 \frac{1}{3}\leq &\frac{D(3,2)}{D(3,8)} \\
 1\leq &\frac{D(3,2)}{D(3,9)} \\
 1\leq &\frac{D(3,2)}{D(3,10)} \\
 &\frac{D(3,4)}{D(3,5)}=2 \\
 1\leq &\frac{D(3,4)}{D(3,6)}\leq 3 \\
 1\leq &\frac{D(3,4)}{D(3,7)}\leq 3 \\
 1\leq &\frac{D(3,4)}{D(3,8)} \\
 3\leq &\frac{D(3,4)}{D(3,9)} \\
 3\leq &\frac{D(3,4)}{D(3,10)} \\
 \frac{1}{2}\leq &\frac{D(3,5)}{D(3,6)}\leq \frac{3}{2} \\ \\ \\
\end{aligned}$
\end{center}
\end{minipage}
&
\begin{minipage}{0.3\linewidth}
\ \ \ $\begin{aligned}
 \frac{1}{2}\leq &\frac{D(3,5)}{D(3,7)}\leq \frac{3}{2} \\
 \frac{1}{2}\leq &\frac{D(3,5)}{D(3,8)} \\
 \frac{3}{2}\leq &\frac{D(3,5)}{D(3,9)} \\
 \frac{3}{2}\leq &\frac{D(3,5)}{D(3,10)} \\
 C_1\leq &\frac{D(3,6)}{D(3,7)}\leq C_2 \\
 1\leq &\frac{D(3,6)}{D(3,8)} \\
 1\leq &\frac{D(3,6)}{D(3,9)} \\
 1\leq &\frac{D(3,6)}{D(3,10)} \\
 1\leq &\frac{D(3,7)}{D(3,8)} \\
 1\leq &\frac{D(3,7)}{D(3,9)} \\
 1\leq &\frac{D(3,7)}{D(3,10)} \\
 \frac{1}{2}\leq &\frac{D(3,9)}{D(3,10)}\leq 1 \\
\end{aligned}$
\end{minipage}
\end{tabular}
\end{center}

\medskip
where
$C_1\approx 0.9002270330$ is the second largest root of
$$6137 x^5-14689 x^4+14429 x^3-9547 x^2+3698 x-100$$
and
$C_2\approx 1.100851119$ is the largest root of the same polynomial.

\newpage

\begin{theorem}
The following bounds involving distances from $X_4$ hold for all triangles.
\end{theorem}

\begin{center}
\begin{tabular}{c|c|c}
\begin{minipage}{0.3\linewidth}
\begin{center}
$\begin{aligned}
&\frac{D(4,1)}{D(4,2)}\leq 1 \\
 &\frac{D(4,1)}{D(4,3)}\leq \frac{2}{3} \\
 &\frac{D(4,1)}{D(4,5)}\leq \frac{4}{3} \\
 1\leq &\frac{D(4,1)}{D(4,6)} \\
 1\leq &\frac{D(4,1)}{D(4,7)}\leq 2 \\
 &\frac{D(4,1)}{D(4,8)}\leq 1 \\
 &\frac{D(4,1)}{D(4,9)}\leq 1 \\
 &\frac{D(4,1)}{D(4,10)}\leq 1 \\
 &\frac{D(4,2)}{D(4,3)}=\frac{2}{3} \\
 &\frac{D(4,2)}{D(4,5)}=\frac{4}{3} \\
 1\leq &\frac{D(4,2)}{D(4,6)} \\
 1\leq &\frac{D(4,2)}{D(4,7)} \\
\end{aligned}$
\end{center}
\end{minipage}
&
\begin{minipage}{0.26\linewidth}
\begin{center}
$\begin{aligned}
 \frac{1}{3}\leq &\frac{D(4,2)}{D(4,8)}\leq 1 \\
 \frac{2}{3}\leq &\frac{D(4,2)}{D(4,9)}\leq 1 \\
 \frac{2}{3}\leq &\frac{D(4,2)}{D(4,10)}\leq 1 \\
 &\frac{D(4,3)}{D(4,5)}=2 \\
 \frac{3}{2}\leq &\frac{D(4,3)}{D(4,6)} \\
 \frac{3}{2}\leq &\frac{D(4,3)}{D(4,7)} \\
 \frac{1}{2}\leq &\frac{D(4,3)}{D(4,8)}\leq \frac{3}{2} \\
 1\leq &\frac{D(4,3)}{D(4,9)}\leq \frac{3}{2} \\
 1\leq &\frac{D(4,3)}{D(4,10)}\leq \frac{3}{2} \\
 \frac{3}{4}\leq &\frac{D(4,5)}{D(4,6)} \\
 \frac{3}{4}\leq &\frac{D(4,5)}{D(4,7)} \\
 \frac{1}{4}\leq &\frac{D(4,5)}{D(4,8)}\leq \frac{3}{4} \\
\end{aligned}$
\end{center}
\end{minipage}
&
\begin{minipage}{0.3\linewidth}
$\begin{aligned}
 \frac{1}{2}\leq &\frac{D(4,5)}{D(4,9)}\leq \frac{3}{4} \\
 \frac{1}{2}\leq &\frac{D(4,5)}{D(4,10)}\leq \frac{3}{4} \\
 &\frac{D(4,6)}{D(4,7)}\leq
   C_3\\
 &\frac{D(4,6)}{D(4,8)}\leq 1 \\
 &\frac{D(4,6)}{D(4,9)}\leq 1 \\
 &\frac{D(4,6)}{D(4,10)}\leq 1 \\
 &\frac{D(4,7)}{D(4,8)}\leq 1 \\
 &\frac{D(4,7)}{D(4,9)}\leq 1 \\
 &\frac{D(4,7)}{D(4,10)}\leq 1 \\
 1\leq &\frac{D(4,8)}{D(4,9)}\leq 2 \\
 1\leq &\frac{D(4,8)}{D(4,10)}\leq 2 \\
 1\leq &\frac{D(4,9)}{D(4,10)}\leq \frac{10}{9} \\
\end{aligned}$
\end{minipage}
\end{tabular}
\end{center}

\medskip
where $C_3\approx1.104068697$ is the positive root of
$8 x^4-36 x^3+113 x^2-69 x-25$.

\newpage

\void{
\begin{minipage}{5in}
\twocolumn
\hspace{0.5in}
$\begin{aligned}
\end{aligned}$$
\hspace{1in}
$$\begin{aligned}
\end{aligned}$
\end{minipage}
\end{theorem}
\hbox{where $C_3\approx1.104068697$ is the positive root of
$8 x^4-36 x^3+113 x^2-69 x-25$.}
\onecolumn
}

\begin{theorem}
The following bounds involving distances from $X_5$ hold for all triangles.
\end{theorem}

\begin{center}
\begin{tabular}{c|c|c}
\begin{minipage}{0.3\linewidth}
\begin{center}
$\begin{aligned}
&\frac{D(5,1)}{D(5,2)}\leq 3 \\
 &\frac{D(5,1)}{D(5,3)}\leq 1 \\
 &\frac{D(5,1)}{D(5,4)}\leq 1 \\
 &\frac{D(5,1)}{D(5,8)}\leq 1 \\
 &\frac{D(5,1)}{D(5,9)}\leq 1 \\
 &\frac{D(5,1)}{D(5,10)}\leq 1 \\
 &\frac{D(5,2)}{D(5,3)}=\frac{1}{3} \\
 &\frac{D(5,2)}{D(5,4)}=\frac{1}{3} \\
 \frac{1}{3}\leq &\frac{D(5,2)}{D(5,6)} \\
 \frac{1}{3}\leq &\frac{D(5,2)}{D(5,7)} \\
 \frac{1}{9}\leq &\frac{D(5,2)}{D(5,8)}\leq 1 \\
\end{aligned}$
\end{center}
\end{minipage}
&
\begin{minipage}{0.26\linewidth}
\begin{center}
$\begin{aligned}
 \frac{1}{3}\leq &\frac{D(5,2)}{D(5,9)}\leq 1 \\
 \frac{1}{3}\leq &\frac{D(5,2)}{D(5,10)}\leq 1 \\
 &\frac{D(5,3)}{D(5,4)}=1 \\
 1\leq &\frac{D(5,3)}{D(5,6)} \\
 1\leq &\frac{D(5,3)}{D(5,7)} \\
 \frac{1}{3}\leq &\frac{D(5,3)}{D(5,8)}\leq 3 \\
 1\leq &\frac{D(5,3)}{D(5,9)}\leq 3 \\
 1\leq &\frac{D(5,3)}{D(5,10)}\leq 3 \\
 1\leq &\frac{D(5,4)}{D(5,6)} \\
 1\leq &\frac{D(5,4)}{D(5,7)} \\
 \frac{1}{3}\leq &\frac{D(5,4)}{D(5,8)}\leq 3 \\
\end{aligned}$
\end{center}
\end{minipage}
&
\begin{minipage}{0.35\linewidth}
\ \ \ $\begin{aligned}
 1\leq &\frac{D(5,4)}{D(5,9)}\leq 3 \\
 1\leq &\frac{D(5,4)}{D(5,10)}\leq 3 \\
 &\frac{D(5,6)}{D(5,8)}\leq 1 \\
 &\frac{D(5,6)}{D(5,9)}\leq C_4 \\
 &\frac{D(5,6)}{D(5,10)}\leq C_5 \\
 &\frac{D(5,7)}{D(5,8)}\leq 1 \\
 &\frac{D(5,7)}{D(5,9)}\leq 1 \\
 &\frac{D(5,7)}{D(5,10)}\leq 1 \\
 1\leq &\frac{D(5,8)}{D(5,9)}\leq 3 \\
 1\leq &\frac{D(5,8)}{D(5,10)}\leq 3 \\
 1\leq &\frac{D(5,9)}{D(5,10)}\leq 7-4 \sqrt{2} \\
\end{aligned}$
\end{minipage}
\end{tabular}
\end{center}

\medskip
where
$C_4\approx 1.053322135$ is the positive root of
$$6137 x^5+5335 x^4+678 x^3-3702 x^2-9479 x-1225$$
and
$C_5\approx 1.194505073$ is the positive root of
$$x^4+2 x^3+22 x^2-30 x-1.$$

\newpage

\begin{theorem}
The following bounds involving distances from $X_6$ hold for all triangles.
\end{theorem}

\begin{center}
\begin{tabular}{c|c|c}
\begin{minipage}{0.3\linewidth}
\begin{center}
$\begin{aligned}
 &\frac{D(6,1)}{D(6,2)}\leq 9-6 \sqrt{2} \\
 &\frac{D(6,1)}{D(6,3)}\leq 2-\sqrt{3} \\
 \frac{1}{2}\leq &\frac{D(6,1)}{D(6,7)} \\
 &\frac{D(6,1)}{D(6,8)}\leq \frac{2 \sqrt{2}-1}{7} \\
 &\frac{D(6,1)}{D(6,9)}\leq \frac{1}{3} \\
 &\frac{D(6,1)}{D(6,10)}\leq \sqrt{2}-1 \\
 &\frac{D(6,2)}{D(6,3)}\leq \frac{2}{3} \\
 1\leq &\frac{D(6,2)}{D(6,7)} \\
\end{aligned}$
\end{center}
\end{minipage}
&
\begin{minipage}{0.33\linewidth}
\begin{center}
$\begin{aligned}
 \frac{1}{3}\leq &\frac{D(6,2)}{D(6,8)}\leq \frac{5+4 \sqrt{2}}{21} \\
 \frac{1}{2}\leq &\frac{D(6,2)}{D(6,9)}\leq \frac{3}{4} \\
 \frac{2}{3}\leq &\frac{D(6,2)}{D(6,10)}\leq \frac{1+\sqrt{2}}{3} \\
 \frac{1}{2}\leq &\frac{D(6,3)}{D(6,4)} \\
 2\leq &\frac{D(6,3)}{D(6,5)} \\
 C_6\leq &\frac{D(6,3)}{D(6,7)} \\
 \frac{1}{2}\leq &\frac{D(6,3)}{D(6,8)} \\
 1\leq &\frac{D(6,3)}{D(6,9)} \\
\end{aligned}$
\end{center}
\end{minipage}
&
\begin{minipage}{0.3\linewidth}
$\begin{aligned}
 1\leq &\frac{D(6,3)}{D(6,10)} \\
 &\frac{D(6,7)}{D(6,8)}\leq \frac{1}{2} \\
 &\frac{D(6,7)}{D(6,9)}\leq \frac{1}{2} \\
 &\frac{D(6,7)}{D(6,10)}\leq \frac{4}{5} \\
 1\leq &\frac{D(6,8)}{D(6,9)}\leq 2 \\
 3-\sqrt{2}\leq &\frac{D(6,8)}{D(6,10)}\leq 2 \\
 1\leq &\frac{D(6,9)}{D(6,10)}\leq \frac{8}{5} \\ \\ \\
\end{aligned}$
\end{minipage}
\end{tabular}
\end{center}

\medskip
where
$C_6\approx 7.8631112181$ is the largest root of
$$6967296 x^{20}+2015974656 x^{18}-160813808784 x^{16}+603818269839 x^{14}-894980577861
   x^{12}$$
$${}+677249814873 x^{10}-274035844587 x^8+60418557684 x^6-6782842860 x^4+290960784
   x^2$$
$+7744$.

\newpage

\begin{theorem}
The following bounds involving distances from $X_7$ hold for all triangles.
\end{theorem}

\begin{center}
\begin{tabular}{c|c|c}
\begin{minipage}{0.4\linewidth}
\begin{center}
$\begin{aligned}
&\frac{D(7,1)}{D(7,2)}\leq \frac{3}{4} \\
 &\frac{D(7,1)}{D(7,3)}\leq \frac{2}{17} \left(5-2 \sqrt{2}\right) \\
 &\frac{D(7,1)}{D(7,4)}\leq 1 \\
 &\frac{D(7,1)}{D(7,8)}\leq \frac{1}{2} \\
 &\frac{D(7,1)}{D(7,9)}\leq \frac{1}{2} \\
 &\frac{D(7,1)}{D(7,10)}\leq \frac{2}{3} \\
 &\frac{D(7,2)}{D(7,3)}\leq \frac{2}{3} \\
 2\leq &\frac{D(7,2)}{D(7,6)} \\
\end{aligned}$
\end{center}
\end{minipage}
&
\begin{minipage}{0.26\linewidth}
\begin{center}
$\begin{aligned}
 \frac{1}{3}\leq &\frac{D(7,2)}{D(7,8)}\leq \frac{2}{3} \\
 &\frac{D(7,2)}{D(7,9)}=\frac{2}{3} \\
 \frac{2}{3}\leq &\frac{D(7,2)}{D(7,10)}\leq \frac{8}{9} \\
 \frac{1}{2}\leq &\frac{D(7,3)}{D(7,4)} \\
 2\leq &\frac{D(7,3)}{D(7,5)} \\
 C_7\leq &\frac{D(7,3)}{D(7,6)} \\
 \frac{1}{2}\leq &\frac{D(7,3)}{D(7,8)} \\
 1\leq &\frac{D(7,3)}{D(7,9)} \\
\end{aligned}$
\end{center}
\end{minipage}
&
\begin{minipage}{0.3\linewidth}
\ \ \ $\begin{aligned}
 1\leq &\frac{D(7,3)}{D(7,10)} \\
 1\leq &\frac{D(7,4)}{D(7,6)} \\
 &\frac{D(7,6)}{D(7,8)}\leq \frac{1}{3} \\
 &\frac{D(7,6)}{D(7,9)}\leq \frac{1}{3} \\
 &\frac{D(7,6)}{D(7,10)}\leq \frac{4}{9} \\
 1\leq &\frac{D(7,8)}{D(7,9)}\leq 2 \\
 \frac{4}{3}\leq &\frac{D(7,8)}{D(7,10)}\leq 2 \\
 1\leq &\frac{D(7,9)}{D(7,10)}\leq \frac{4}{3} \\
\end{aligned}$
\end{minipage}
\end{tabular}
\end{center}

\medskip
where
$C_7\approx 7.9776615835$ is the largest root of\\
{\raggedright$833089536x^{28}+220028016384 x^{26}-19474287964848 x^{24}+139707882692901
   x^{22}$ \hspace*{12pt}${}-410390834384412 x^{20}+732430210466916 x^{18}-892396597211316
   x^{16}$ \hspace*{12pt}${}+782711166381062 x^{14}-492062343977916 x^{12}+216425700787620
   x^{10}$\\ \hspace*{12pt}${}-65960002546284 x^8+14226627485565 x^6-2259294716376 x^4+253570773456
   x^2$ \hspace*{12pt}${}-14637417984$.}

\newpage

\begin{theorem}
The following bounds involving distances from $X_8$ hold for all triangles.
\end{theorem}

\begin{center}
\begin{tabular}{c|c}
\begin{minipage}{0.45\linewidth}
\begin{center}
$\begin{aligned}
 &\frac{D(8,1)}{D(8,2)}=\frac{3}{2} \\
 &\frac{D(8,1)}{D(8,4)}\leq 1 \\
 &\frac{D(8,1)}{D(8,5)}\leq \frac{4}{3} \\
 \frac{2}{7} \left(4-\sqrt{2}\right)\leq &\frac{D(8,1)}{D(8,6)}\leq 1 \\
 \frac{1}{2}\leq &\frac{D(8,1)}{D(8,7)}\leq 1 \\
 2\leq &\frac{D(8,1)}{D(8,9)} \\
 &\frac{D(8,1)}{D(8,10)}=2 \\
 &\frac{D(8,2)}{D(8,4)}\leq \frac{2}{3} \\
 &\frac{D(8,2)}{D(8,5)}\leq \frac{8}{9} \\
 \frac{4}{21} \left(4-\sqrt{2}\right)\leq &\frac{D(8,2)}{D(8,6)}\leq \frac{2}{3} \\
 \frac{1}{3}\leq &\frac{D(8,2)}{D(8,7)}\leq \frac{2}{3} \\
 \frac{4}{3}\leq &\frac{D(8,2)}{D(8,9)} \\
 &\frac{D(8,2)}{D(8,10)}=\frac{4}{3} \\
 &\frac{D(8,3)}{D(8,4)}\leq \frac{1}{2} \\
 &\frac{D(8,3)}{D(8,5)}\leq 2 \\
\end{aligned}$
\end{center}
\end{minipage}
&
\begin{minipage}{0.5\linewidth}
\begin{center}
$\begin{aligned}
 \frac{4}{3}\leq &\frac{D(8,4)}{D(8,5)}\leq 4 \\
 1\leq &\frac{D(8,4)}{D(8,6)} \\
 1\leq &\frac{D(8,4)}{D(8,7)} \\
 2\leq &\frac{D(8,4)}{D(8,9)} \\
 2\leq &\frac{D(8,4)}{D(8,10)} \\
 C_8\leq
   &\frac{D(8,5)}{D(8,6)} \\
 \frac{1}{8} \left(3+2 \sqrt{2}\right)\leq &\frac{D(8,5)}{D(8,7)} \\
 \frac{3}{2}\leq &\frac{D(8,5)}{D(8,9)} \\
 \frac{3}{2}\leq &\frac{D(8,5)}{D(8,10)} \\
 \frac{2}{3}\leq &\frac{D(8,6)}{D(8,7)}\leq \frac{7}{6} \\
 2\leq &\frac{D(8,6)}{D(8,9)} \\
 2\leq &\frac{D(8,6)}{D(8,10)}\leq 2+\frac{1}{\sqrt{2}} \\
 2\leq &\frac{D(8,7)}{D(8,9)} \\
 2\leq &\frac{D(8,7)}{D(8,10)}\leq 4 \\
 &\frac{D(8,9)}{D(8,10)}\leq 1 \\
\end{aligned}$
\end{center}
\end{minipage}
\end{tabular}
\end{center}

\medskip
where
$C_8\approx 0.6817039304$ is the smallest positive root of
$$896 x^4-2184 x^3+1924 x^2-758 x+121.$$

\newpage

\begin{theorem}
The following bounds involving distances from $X_{9}$ hold for all triangles.
\end{theorem}

\begin{center}
\begin{tabular}{c|c}
\begin{minipage}{0.4\linewidth}
\begin{center}
$\begin{aligned}
 \frac{3}{2}\leq &\frac{D(9,1)}{D(9,2)}\leq 3 \\
 &\frac{D(9,1)}{D(9,4)}\leq 1 \\
 &\frac{D(9,1)}{D(9,5)}\leq 2 \\
 \frac{2}{3}\leq &\frac{D(9,1)}{D(9,6)}\leq 1 \\
 \frac{1}{2}\leq &\frac{D(9,1)}{D(9,7)}\leq 1 \\
 1\leq &\frac{D(9,1)}{D(9,8)} \\
 2\leq &\frac{D(9,1)}{D(9,10)} \\
 &\frac{D(9,2)}{D(9,4)}\leq \frac{1}{3} \\
 &\frac{D(9,2)}{D(9,5)}\leq \frac{2}{3} \\
 \frac{1}{4}\leq &\frac{D(9,2)}{D(9,6)}\leq \frac{1}{2} \\
 &\frac{D(9,2)}{D(9,7)}=\frac{1}{3} \\
 \frac{1}{3}\leq &\frac{D(9,2)}{D(9,8)} \\
 \frac{4}{3}\leq &\frac{D(9,2)}{D(9,10)} \\
 &\frac{D(9,3)}{D(9,4)}\leq \frac{1}{2} \\
 &\frac{D(9,3)}{D(9,5)}\leq 2 \\
\end{aligned}$
\end{center}
\end{minipage}
&
\begin{minipage}{0.4\linewidth}
\begin{center}
$\begin{aligned}
 1\leq &\frac{D(9,3)}{D(9,10)} \\
 2\leq &\frac{D(9,4)}{D(9,5)}\leq 4 \\
 1\leq &\frac{D(9,4)}{D(9,6)} \\
 1\leq &\frac{D(9,4)}{D(9,7)} \\
 1\leq &\frac{D(9,4)}{D(9,8)} \\
 10\leq &\frac{D(9,4)}{D(9,10)} \\
 C_9\leq &\frac{D(9,5)}{D(9,6)} \\
 \frac{1}{2}\leq &\frac{D(9,5)}{D(9,7)} \\
 \frac{1}{2}\leq &\frac{D(9,5)}{D(9,8)} \\
 \frac{5}{2}+\sqrt{2}\leq &\frac{D(9,5)}{D(9,10)} \\
 \frac{2}{3}\leq &\frac{D(9,6)}{D(9,7)}\leq \frac{4}{3} \\
 1\leq &\frac{D(9,6)}{D(9,8)} \\
 \frac{8}{3}\leq &\frac{D(9,6)}{D(9,10)} \\
 1\leq &\frac{D(9,7)}{D(9,8)} \\
 4\leq &\frac{D(9,7)}{D(9,10)} \\
\end{aligned}$
\end{center}
\end{minipage}
\end{tabular}
\end{center}

\medskip
where
$C_9\approx 0.4870156430$ is the smallest positive root of
$$3072 x^5+9304 x^4-35096 x^3+40708 x^2-25350 x+6137.$$

\newpage

\begin{theorem}
The following bounds involving distances from $X_{10}$ hold for all triangles.
\end{theorem}

\begin{center}
\begin{tabular}{c|c}
\begin{minipage}{0.45\linewidth}
\begin{center}
$\begin{aligned}
 &\frac{D(10,1)}{D(10,2)}=3 \\
 &\frac{D(10,1)}{D(10,4)}\leq 1 \\
 &\frac{D(10,1)}{D(10,5)}\leq 2 \\
 2-\sqrt{2}\leq &\frac{D(10,1)}{D(10,6)}\leq 1 \\
 \frac{1}{3}\leq &\frac{D(10,1)}{D(10,7)}\leq 1 \\
 &\frac{D(10,1)}{D(10,8)}=1 \\
 1\leq &\frac{D(10,1)}{D(10,9)} \\
 &\frac{D(10,2)}{D(10,4)}\leq \frac{1}{3} \\
 &\frac{D(10,2)}{D(10,5)}\leq \frac{2}{3} \\
 \frac{1}{3} \left(2-\sqrt{2}\right)\leq &\frac{D(10,2)}{D(10,6)}\leq \frac{1}{3} \\
 \frac{1}{9}\leq &\frac{D(10,2)}{D(10,7)}\leq \frac{1}{3} \\
 &\frac{D(10,2)}{D(10,8)}=\frac{1}{3} \\
 \frac{1}{3}\leq &\frac{D(10,2)}{D(10,9)} \\
 &\frac{D(10,3)}{D(10,4)}\leq \frac{1}{2} \\
 &\frac{D(10,3)}{D(10,5)}\leq 2 \\
 2\leq &\frac{D(10,3)}{D(10,9)} \\
\end{aligned}$
\end{center}
\end{minipage}
&
\begin{minipage}{0.45\linewidth}
\begin{center}
$\begin{aligned}
 2\leq &\frac{D(10,4)}{D(10,5)}\leq 4 \\
 1\leq &\frac{D(10,4)}{D(10,6)} \\
 1\leq &\frac{D(10,4)}{D(10,7)} \\
 1\leq &\frac{D(10,4)}{D(10,8)} \\
 9\leq &\frac{D(10,4)}{D(10,9)} \\
 C_{10}\leq
   &\frac{D(10,5)}{D(10,6)} \\
 \frac{1}{2}\leq &\frac{D(10,5)}{D(10,7)} \\
 \frac{1}{2}\leq &\frac{D(10,5)}{D(10,8)} \\
 \frac{3}{2}+\sqrt{2}\leq &\frac{D(10,5)}{D(10,9)} \\
 \frac{5}{9}\leq &\frac{D(10,6)}{D(10,7)}\leq \frac{4}{3} \\
 1\leq &\frac{D(10,6)}{D(10,8)}\leq 1+\frac{1}{\sqrt{2}} \\
 \frac{5}{3}\leq &\frac{D(10,6)}{D(10,9)} \\
 1\leq &\frac{D(10,7)}{D(10,8)}\leq 3 \\
 3\leq &\frac{D(10,7)}{D(10,9)} \\
 1\leq &\frac{D(10,8)}{D(10,9)} \\ \\ \\
\end{aligned}$
\end{center}
\end{minipage}
\end{tabular}
\end{center}

\medskip
where
$C_{10}\approx 0.4870156430$ is the smallest positive root of
$$50 x^4-72 x^3+22 x^2-2 x+1.$$

\newpage


\end{document}